\documentclass[10pt,a4paper]{article}
\usepackage{amsmath}
\usepackage{amsfonts}
\usepackage{amssymb}
\usepackage{graphicx}
\usepackage{amsthm}
\usepackage[all]{xy}
\usepackage{enumitem}

\usepackage[subnum]{cases}
\usepackage[colorlinks=true,linkcolor=blue,citecolor=blue]{hyperref}

\newcommand{\R}{\mathbb{R}}
\newcommand{\dd}{\mathrm{d}}
\newcommand{\Lie}{\mathcal{L}}
\newcommand{\id}{\mathrm{id}}
\newcommand{\Gb}{\mathcal{G}}
\newcommand{\Ab}{\mathcal{A}}

\newcommand{\Dd}{\mathrm{D}}

\newcommand{\bX}{\mathfrak{X}_b}

\newcounter{oftheorem}[subsection]
\newenvironment{mytheorem}[1]%
{\begin{trivlist}
     \renewcommand{\theoftheorem}{#1~\thesection.\arabic{oftheorem}}
     \refstepcounter{oftheorem}
     \item[\hspace{\labelsep}\bf\thesection.\arabic{oftheorem} #1]}%
{\end{trivlist}}
\newenvironment{definition}{\begin{mytheorem}{Definition}\it}{\end{mytheorem}}

\newenvironment{proposition}{\begin{mytheorem}{Proposition}\it}{\end{mytheorem}}
\newenvironment{theorem}{\begin{mytheorem}{Theorem}\it}{\end{mytheorem}}
\newenvironment{corollary}{\begin{mytheorem}{Corollary}\it}{\end{mytheorem}}
\newenvironment{remark}{\begin{mytheorem}{Remark}}{\end{mytheorem}}
\newenvironment{lemma}{\begin{mytheorem}{Lemma}}{\end{mytheorem}}

\author{Stavros Anastassiou\\Department of Mathematics\\
 University of Western Macedonia\\Kastoria, Greece\\sanastassiou@gmail.com}

\date{}

\title{Local study of generic dynamical systems on the b-symplectic plane}

\begin{document}

\maketitle

\begin{abstract}
Let $\Gb$ be the group of germs of planar diffeomorphisms, fixing the origin, which preserve the standard b-symplectic form. The group acts on itself, under the usual left-right action and we identify, in the generic case, a simple representative of each orbit, depending on  a non-vanishing, univariate, function. We then study the action of $\Gb$ on b-symplectic vector fields. For those fields that do not vanish at the origin, we establish a flow-box type theorem, while, for those that do vanish there, we show that the study of their conjugacy classes reduces to a classification of functions. We present this classification, when the origin is a regular  point for these functions, while if it is a critical point, we present the orbits of their 2-jets. Finally, we address the problem of whether an element of $\Gb$ can be realised as the time-1 map of a b-symplectic vector field.
\end{abstract}


\section{Introduction}\label{sec:intro}
Let $M$ be a smooth (i.e. $C^{\infty}$) manifold of dimension $2n$ and $Z\subset M$ a closed hypersurface. A \emph{$b$-vector field} is a vector field on $M$ which is tangent to $Z$. They form the $C^{\infty}(M)$-module $\bX M$. By the Serre--Swan theorem (\cite{Nestruev}) there exists a vector bundle $^bT M\to M$ whose sections are exactly $\bX(M)$, the so-called \emph{$b$-tangent bundle}. 

The \emph{$b$-cotangent bundle} is $^bT^*M=(^bT M)^*$, while the space of \emph{b-k-forms} is $^b \Omega ^k(M)=\Gamma(\Lambda ^k\,^{b}T^* M)$.

\begin{definition}\label{def:bsymp}
A \emph{$b$-symplectic form} on $(M^{2n},Z)$ is a closed b-2-form $\omega\in ^b \Omega ^2(M)$, which is nondegenerate as a section of $\Lambda^2\,^bT^* M$ at every point of $M$ (equivalently, $\omega^n$ is a nowhere-vanishing section of $\Lambda^{2n}\ ^bT^* M$).
\end{definition}

On $M \setminus Z$, $\omega$ restricts to an ordinary symplectic form with a log-symplectic singularity (or simple pole) along $Z$.

The Darboux theorem for $b$-symplectic forms (\cite{GMP}) states that, if $p\in Z$, there exist coordinates $(x_1,y_1,\dots,x_{n-1},y_{n-1},z,t)$, centred at $p$, such that $Z=\{z=0\}$ and:
\[
\omega=\sum_{i=1}^{n-1}\dd x_i\wedge\dd y_i+\frac{\dd z}{z}\wedge\dd t .
\]

In this article, we focus on the two-dimensional case and study the problem of classifying diffeomorphisms and vector fields which preserve this b-symplectic form locally, in a neighbourhood of a point belonging to $Z$. The diffeomorphisms which preserve $\omega$ form a group $\Gb$ and we study their general form in section 2, where we also distinguish a subset of $\Gb$, the mappings of which are quite simple to study. In Section 3, after introducing the natural action of $\Gb$ on itself, we prove that a generic member of $\Gb$ is conjugate to one of these simple mappings, thereby providing explicit local models for almost all elements of $\Gb$, using classical results from invariant manifold theory. In section 4 we study the analogous problem for vector fields preserving the b-symplectic form. We first establish a flow-box theorem for non-vanishing vector fields. For vector fields vanishing at the point of interest, we show that their classification descents to a classification of functions. We proceed to present this classification, in the non-critical case. Although functional moduli appear, we are able to construct simple local models for these functions and their corresponding vector fields. At critical points, we describe the orbits of their 2-jets and compare our results with a classification of Arnold, regarding singularities of functions on manifolds with boundary. In section 5, we solve the embedding problem: we explicitly determine necessary and sufficient conditions for a $b$-symplectic diffeomorphism to be the time-1 map of a $b$-symplectic vector field. Finally, Section 6 contains concluding remarks and discusses the extension of our results to higher dimensions.

For the general study of b-symplectic manifolds the interested reader should consult \cite{BDMOP,DKM,GMP,KMS}. Article \cite{Kou} is also related with a classification of functions similar to our own, although in the complex-analytic setting, while \cite{Kou15} addresses similar problems, when a symplectic structure vanishes along a line of the plane.

\section{Symmetries of the b-symplectic plane}\label{sec:symm}

Equip the plane $\R^2$ with coordinates $(z,t)$, let $Z=\{z=0\}$ and, on $\R^2\setminus Z$, define the form:
\[
\omega=\frac{1}{z}\dd z\wedge\dd t,
\]
which is the standard $b$-symplectic form.

\begin{definition}\label{def:G}
We denote by $\Gb$ the group of germs at the origin of diffeomorphisms $\varphi:U\to V$, where $U,V$ are open neighbourhoods of $(0,0)\in\R^2$, such that:
\begin{enumerate}[label=(\roman*),nosep]
\item $\varphi(0,0)=(0,0)$;
\item $\varphi^*\omega=\omega$ whenever both sides are defined, that is, in the open and dense set $O_\varphi:=U\setminus\bigl(Z\cup\varphi^{-1}(Z)\bigr)$.
\end{enumerate}
We shall call the members of $\Gb$ ``b-symplectomorphisms''.
\end{definition}

Although we do not explicitly assume that b-symplectomorphisms preserve the line $Z$, this follows from the definition above.

\begin{lemma}\label{lem:Zinv}
Let $U,V$ be open neighbourhoods of the origin and $\varphi :U\to V$ a b-symplectomorphism.  Then $\varphi(U\cap Z)=V\cap Z$.
\end{lemma}
\begin{proof}
Write $\varphi=(F,T)$. On $O_\varphi$ we have $F\neq0$, because $\varphi(O_\varphi)\cap Z=\emptyset$ and we can compute that:
\[
\varphi^*\omega=\frac{\dd F\wedge\dd T}{F}=\frac{\det D\varphi}{F}\,\dd z\wedge\dd t .
\]
But $\varphi^*\omega=\omega$ on $O_\varphi$ and so we get the equation:
\begin{equation}\label{eq:pointwise}
F(z,t)=z\,\det D\varphi(z,t).
\end{equation}

Now, to prove that $\varphi(U\cap Z)\subset Z$, let $p=(0,t_0)\in U\cap Z$ and suppose $F(p)\neq0$. By continuity, $F\neq0$ near $p$ and so the points $(z,t_0)$ (with $z\neq0$ small enough) lie in $O_\varphi$. Due to the relation \eqref{eq:pointwise}, $|F(z,t_0)|\to0$ as $z\to0$, which give us $F(p)=0$ and we arrive at a contradiction.

$\varphi^{-1}$ is also a b-symplectomorphism and, using the same reasoning, $\varphi^{-1}(V\cap Z)\subset U\cap Z$, i.e.\ $V\cap Z\subset\varphi(U\cap Z)$; therefore $\varphi(U\cap Z)=V\cap Z$.
\end{proof}

\begin{remark}
In particular, the lemma above ensures that $\varphi^{-1}(Z)\cap U=Z\cap U$, $O_\varphi=U\setminus Z$ and $\varphi^*\omega=\omega$ on $U\setminus Z$.
\end{remark}

Although we used equation (\ref{eq:pointwise}) on $O_{\varphi}$, it holds even on $Z$, since both of its sides are continuous functions and $O_{\varphi}$ is dense in $U$. This equation will be needed in the future, and we shall refer to it as ``the Jacobian equation".

Let us now obtain the general form of an element of $\Gb$.

\begin{lemma}\label{lem:structure}
Let $\varphi:U\to V$ represent an element of $\Gb$, so that $U,V\subset\R^2$ are open neighbourhoods of $(0,0)$ and $\varphi(0,0)=(0,0)$. Choose $\delta,\varepsilon>0$ such that the rectangle:
\[
R=R(\delta,\varepsilon):=(-\delta,\delta)\times(-\varepsilon,\varepsilon)
\]
is contained in $U$, and write $\varphi|_R=(F,T)$, where $F,T$ smooth. There exist smooth functions $u,S,\tilde u:R\to\R$ such that:
\begin{enumerate}[label=(\roman*)]
\item $\varphi(z,t)=\bigl(z\,u(z,t),\ t+z\,S(z,t)\bigr)$, with $u(z,t)=u_0(t)+z\,\tilde u(z,t)$ and $u_0(t):=u(0,t)$.

\item $u(z,t)\neq0$ for every $(z,t)\in R$; in particular $u_0(t)\neq0$ for every $t\in(-\varepsilon,\varepsilon)$.
\end{enumerate}
\end{lemma}

\begin{proof}
By \ref{lem:Zinv}, it is $\varphi(U\cap Z)=V\cap Z$. Hence:
\begin{itemize}[nosep]
\item[(a)] $F(0,t)=0$ for all $|t|<\varepsilon$;
\item[(b)] $F(z,t)\neq0$ for $z\neq0$;
\item[(c)] $T(0,0)=0$, since $\varphi(0,0)=(0,0)$.
\end{itemize}
We now prove that $F=zu$. If $(z,t)\in R$ and $s\in[0,1]$, the point $(sz,t)$ lies in $R$. Since $F(0,t)=0$, we have:
\[
F(z,t)=\int_0^1\frac{\dd}{\dd s}F(sz,t)\,\dd s=z\,u(z,t),\qquad u(z,t):=\int_0^1F_z(sz,t)\,\dd s,
\]
where the function $u:R\to\R$ is smooth, due to Hadamard's lemma. By using the point (b) above, $F=zu\neq0$ for $z\neq0$, hence $u(z,t)\neq0$ for $z\neq0$.

We can now compute that, at the point $(0,t)\in R$, the determinant of $D\varphi(0,t)$ equals $u_0(t)\cdot T_t(0,t)$. It should be nonzero, allowing us to deduce that $u_0(t)\neq0$; together with the fact that $u\neq0$, for $z\neq0$, this proves the second assertion.

Also, the Jacobian equation, for $z=0$, reads as $u_0(t)\,T_t(0,t)=u_0(t)$, and since $u_0(t)\ne0$ we get $T_t(0,t)=1$, for all $|t|<\varepsilon$. Since $T(0,0)=0$:
\[
T(0,t)=T(0,0)+\int_0^tT_t(0,s)\,\dd s=t.
\]

Finally, the smooth function $T-t$ vanishes on $R\cap Z$ and we can write:
\[
T(z,t)-t=z\,S(z,t),\qquad S(z,t):=\int_0^1T_z(sz,t)\,\dd s,
\]
with $S:R\to\R$ smooth. The function $u(z,t)-u_0(t)$ is smooth on $R$ and vanishes on $R\cap Z$ and thus:
\[
u(z,t)=u_0(t)+z\,\tilde u(z,t),\qquad\tilde u(z,t):=\int_0^1u_z(sz,t)\,\dd s,
\]
with $\tilde u:R\to\R$ smooth and this concludes the proof.
\end{proof}

The lemma above shows that the origin is not the only fixed point of $\varphi$; the whole line $Z$ consists of fixed points and none of them is hyperbolic (but note that fixed points that do not lie on $Z$ may be hyperbolic).

We now turn our attention to some members of $\Gb$, that will be of great importance in what follows.

Let $\Ab$ be the abelian, multiplicative, group of germs at $0\in\R$ of smooth, nowhere vanishing, univariate functions. For $a\in\Ab$, define $L_a(z,t)=(a(t)z,\,t)$. The following lemma is easy to obtain.

\begin{lemma}\label{lem:L}
If $a\in\Ab$, $L_a^*\omega=\omega$ on $\{z\ne0\}$ and thus $L_a\in\Gb$. Moreover, $L_a\circ L_b=L_{ab}$ and $L_a^{-1}=L_{1/a}$.
\end{lemma}

As we saw, every $\varphi\in\Gb$ is of the form $\varphi=(zu,\,t+zS)$, with $u_0(t)=u(0,t)\neq0$, thus $u_0\in \Ab$. This function $u_0$ corresponding to $\varphi$ is of special importance and this motivates the introduction of the following:

\begin{definition}\label{def:rho}
Let $\varphi\in \Gb$ and $u_0\in \Ab$ be the univariate function corresponding to it. We shall call $\varphi$ ``normally hyperbolic'' if $|u_0(0)|\ne1$, ``normally attracting'' if $|u_0(0)|<1$ and ``normally repelling'' if $|u_0(0)|>1$. We shall also write $u_0=\rho(\varphi)$.
\end{definition}

\begin{lemma}\label{lem:rho}
$\rho:\Gb\to\Ab$ is a group epimorphism. The kernel $\Gb_1:=\ker\rho$ consists of those elements $\varphi\in\Gb$, with $u(0,t)\equiv1$.
\end{lemma}
\begin{proof}
Let $\varphi,\psi\in\Gb$ and use \ref{lem:structure} to write them as $\varphi(z,t)=\bigl(z\,u^\varphi(z,t),\ t+z\,S^\varphi(z,t)\bigr)$ and $\psi(z,t)=\bigl(z\,u^\psi(z,t),\ t+z\,S^\psi(z,t)\bigr)$.

The composition $\varphi\circ\psi$ has as first component the term $z\,u$, where $u(z,t)=u^\psi(z,t)\;u^\varphi\bigl(\psi(z,t)\bigr)$. Evaluating at the point $(0,t)$ and since, by \ref{lem:structure}, it is $\psi(0,t)=(0,t)$, we get:
\[
u_0(t):=u(0,t)=u^\psi(0,t)\,u^\varphi(0,t)=u_0^\psi(t)\,u_0^\varphi(t).
\]
We therefore have $\rho(\varphi\circ\psi)=\rho(\varphi)\,\rho(\psi)$.

For the identity b-symplectomorphism the function $u$ is equal to $1$, so $\rho(\id)=1$ and $\rho(\varphi)\,\rho(\varphi^{-1})=\rho(\varphi\circ\varphi^{-1})=\rho(\id)=1$.

To prove that $\rho$ is surjective, note that if $a\in\Ab$, $L_a(z,t)=(a(t)z,t)\in\Gb$, and $\rho(L_a)=a$.

The description of the kernel is immediate from the definition of $\rho$.
\end{proof}

In the next section, we classify b-symplectomorphisms up to the natural action on them.


\section{Local study of generic b-symplectomorphisms}\label{sec:hyp}

The natural action of the group $\Gb$ on itself is given by $(\Psi,\varphi)\mapsto \Psi \cdot \varphi :=\Psi \circ \varphi \circ \Psi^{-1}$, where $\Psi,\varphi \in \Gb$. We want to describe the conjugacy class of $\varphi$, under this action. 

Generically, for the function $u_0=\rho(\varphi)$ the relation $|u_0(0)|\neq 1$ will hold and in this case the study of $\varphi$ can be reduced to the study of $L_{u_{0}}=(u_0z,t)$.

More precisely, in this section, we shall prove the following:

\begin{theorem}\label{thm:A}
Let us suppose that the germ $\varphi\in\Gb$ is normally hyperbolic and let $u_0=\rho(\varphi)$. There exists a  $\Psi\in\Gb_1$, such that $\Psi\circ\varphi\circ\Psi^{-1}=L_{u_0}$.
\end{theorem}

Let $|u_0(0)|>1$. For the germ $\varphi^{-1}\in\Gb$ we have that $\rho(\varphi^{-1})=1/u_0$ (\ref{lem:rho}), so $|\rho(\varphi^{-1})(0)|<1$. Also, if the b-symplectomorphism $\Psi\in\Gb_1$ satisfies the relation $\Psi\varphi^{-1}\Psi^{-1}=L_{1/u_0}$, then $\Psi\varphi\Psi^{-1}=L_{1/u_0}^{-1}=L_{u_0}$. It therefore suffices to prove \ref{thm:A} just in the case where $\varphi$ is normally attracting.

This is why, in this section, we assume that for the b-symplectomorphism $\varphi$ it is $|u_0(0)|<1$, while we write $\varphi=(zu,\ t+zS)$, with $u_0(t)=u(0,t)\neq 0$, in the rectangle $R=R(\delta,\varepsilon)$, as dictated by \ref{lem:structure}.

The sought-after conjugacy $\Psi$, mentioned in the theorem above, will of course be of form $\Psi=(z\pi_1,\pi_2)$ (since it is a $b$-symplectomorphism), for suitably chosen functions $\pi_1$ and $\pi_2$. Let us first define $\pi_2$.
\subsection{\texorpdfstring{The function $\pi_2$}{The function pi2}}

Comparing the second components of the conjugacy relation $\Psi\circ\varphi=L_{u_0}\circ\Psi$, we get $\pi_2\circ\varphi=\pi_2$ and, since $\Psi\in\Gb$, we must also have $\pi_2(0,t)=t$ (recall \ref{lem:structure}). Thus, the function $\pi_2$ has to be both constant on the orbits of $\varphi$ (which converge to points of $Z$) and equal to the identity when restricted on $Z$.

To continue, we need some facts from the theory of normally hyperbolic invariant manifolds. Let $X$ be a smooth manifold, $\Lambda\subset X$ a smooth submanifold and $\Phi$ a smooth map defined near $\Lambda$, with $\Phi(p)=p$ for all $p\in\Lambda$. It follows that, $\forall p\in \Lambda$, $d_p\Phi$ is the identity on $T_p\Lambda$ and it induces the following linear map:
\[
N_p\Phi:\ T_pX/T_p\Lambda\longrightarrow T_pX/T_p\Lambda .
\]
When, as in our case, $\Lambda$ is a hypersurface, $T_pX/T_p \Lambda$ is one-dimensional and $N_p\Phi$ is nothing more than a multiplication by a real number, the so-called ``normal multiplier of $\Phi$ at $p$''.

\begin{theorem}\label{thm:phase}
Let $X$ be a smooth manifold, $\Lambda\subset X$ a compact, smooth, submanifold without boundary, $O\subset X$ an open
neighbourhood of $\Lambda$, and $\Phi:O\to X$ a smooth diffeomorphism onto an open subset of $X$ such that, for every $p\in\Lambda$:
\begin{enumerate}[label=(\alph*)]
\item $\Phi(p)=p$,
\item all eigenvalues of $N_p\Phi$ have modulus smaller than $1$.
\end{enumerate}
Then, for every integer $k\ge1$, there exists an open neighbourhood $W_k\subset O$ of $\Lambda$, with $\Phi(W_k)\subset W_k$, such that, for every $x\in W_k$, the limit $\pi(x):=\lim_{n\to\infty}\Phi^n(x)$ exists and belongs to $\Lambda$, while the map $\pi:W_k\to\Lambda$ is of class $C^k$.
\end{theorem}

\begin{remark}
This is the existence and smoothness theorem for the stable foliation of a normally attracting invariant manifold consisting entirely of fixed points, which in our case holds for every $k\geq 1$, since the eigenvalue in the direction tangent to the invariant manifold is equal to 1 and the other eigenvalue is strictly less than 1. The classic textbook reference for normally hyperbolic invariant manifolds is \cite{HPS}, but see also \cite{Fen}. The paper \cite{Tak71} deals with $C^\infty$ normal forms near partially hyperbolic fixed points, while the book \cite{Aul} covers the case in which the manifold of interest consists of fixed points, as in our case here. 
\end{remark}

Actually, one can prove that the function $\pi$, given by \ref{thm:phase}, is smooth.
\begin{corollary}\label{cor:Cinfty}
The map $\pi:W_1\to\Lambda$, presented in \ref{thm:phase}, is $C^\infty$, and $\pi\circ\Phi=\pi$ on $W_1$.
\end{corollary}
\begin{proof}
Since $\Phi(W_1)\subset W_1$, it is:
\[
\pi(\Phi(x))=\lim_{n\to\infty}\Phi^{n+1}(x)=\pi(x),
\]
for $x\in W_1$. Fix $k\ge1$ and $x_0\in W_1$.
The point $\pi(x_0)=\lim\Phi^n(x_0)$ lies in $\Lambda\subset W_k$ and $W_k$ is open, so there exists $n_0\in \mathbb{N}$, such that $\Phi^{n_0}(x_0)\in W_k$. Due to the fact that $\Phi^{n_0}$ is continuous, there exists a neighbourhood $V\subset W_1$ of $x_0$, such that $\Phi^{n_0}(V)\subset W_k\cap W_1$. On $V$ we have $\pi=\pi\circ\Phi^{n_0}$, and $\pi|_{W_k\cap W_1}$ is $C^k$ (on $W_k\cap W_1$ both definitions of $\pi$ give the same limit), so $\pi|_V$ is $C^k$. But $k$ is arbitrary, so we can conclude that $\pi$ is $C^\infty$ near $x_0$.
\end{proof}

\ref{thm:phase} is stated for a compact $\Lambda$ and we want to apply it, locally around zero, to the non-compact line $Z$ and the germ $\varphi$. To achieve that, we modify $\varphi$ away from the origin and study it in a neighbourhood of a circle, instead of $Z$.

\begin{proposition}\label{prop:A}
Let $\varphi\in\Gb$ be normally attracting. Then there exists a smooth function $\pi$, defined on a neighbourhood of $0\in\R^2$, such that $\pi(0,t)=t$ and $\pi\circ\varphi=\pi$.
\end{proposition}
\begin{proof}
The function $u_0=\rho(\phi)$ is continuous, does not vanish and satisfies $|u_0(0)|<1$. Thus, we can choose an $\varepsilon_0>0$, with $4\varepsilon_0<\varepsilon$, such that $u_0$ has constant sign and satisfies the inequality $|u_0|<1$ on $[-2\varepsilon_0,2\varepsilon_0]$. Set $\ell=8\varepsilon_0$ and $X=(-\delta_1,\delta_1)\times (\mathbb{R}/ \ell \mathbb{Z})$, with coordinates $(z,t)$, $t\in\R / \ell \mathbb{Z}$, where the constant $\delta_1\in(0,\delta]$ will be chosen in a while. We represent classes in $\R /\ell \mathbb{Z}$ by numbers $t\in(-4\varepsilon_0,4\varepsilon_0]$. Let also $\Lambda=\{z=0\}$.

Now, if $\vartheta:\R/\ell\mathbb{Z}\to[0,1]$ is a smooth function, with $\vartheta=1$ on $[-\varepsilon_0,\varepsilon_0]$ and support contained in $(-2\varepsilon_0,2\varepsilon_0)$, define:
\[
\hat{u}(z,t)=\vartheta(t)\,u(z,t)+\bigl(1-\vartheta(t)\bigr)\,u_0(0),\qquad \hat{S}(z,t)=\vartheta(t)\,S(z,t).
\]
The constant $\delta_1$ above was chosen so as $|\hat{u}|<1$ on $X$.

We shall show that the mapping $\Phi$, defined as:
\[
\Phi:X\to X,\qquad\Phi(z,t)=\bigl(z\,\hat{u}(z,t),\ \ t+z\,\hat{S}(z,t)\ \bmod\ell\bigr)
\]
satisfies the hypotheses of \ref{thm:phase}.

\begin{enumerate}[label=(\arabic*)]
\item Obviously, $\Phi(0,t)=(0,t)$ for all $t$.
\item The Jacobian matrix of $\Phi$, at the point $p=(0,t)\in\Lambda$, equals:
\[
D\Phi(p)=\begin{pmatrix}\hat{u}_0(t)&0\\ \hat{S}(0,t)&1\end{pmatrix}
\]
where $\hat{u}_0(t):=\hat{u}(0,t)=\vartheta(t)u_0(t)+(1-\vartheta(t))u_0(0)$. For $|t|<2\varepsilon_0$ this is a convex combination of two numbers of the same sign and of modulus less than $1$, because $u_0$ has constant sign on $[-2\varepsilon_0,2\varepsilon_0]$; for $|t|\ge2\varepsilon_0$ it equals $u_0(0)$. Thus $0<|\hat{u}_0|<1$ everywhere on $\Lambda$. Hence $\Phi$ is a local diffeomorphism at every point of $\Lambda$, with normal multiplier at $p$ equal to $\hat{u}_0(t)$. Since $\Lambda$ is compact and $\Phi|_\Lambda=\id$ is injective, $\Phi$ is injective on a neighbourhood of $\Lambda$, hence a diffeomorphism of an open neighbourhood $O$ of $\Lambda$ onto its image. This proves point (b) of \ref{thm:phase}.
\end{enumerate}

Therefore, \ref{thm:phase} and \ref{cor:Cinfty} apply to $\Phi:O\to \Phi(O)$, to give us an open neighbourhood $W\subset O$ of $\Lambda$ with
$\Phi(W)\subset W$ and a smooth map $\tilde\pi:W\to\Lambda=\R/\ell\mathbb{Z}$, such that:
\[
\tilde\pi\circ\Phi=\tilde\pi\ \text{ on }W,\qquad\tilde\pi(0,t)=t\ \bmod\ell\ \text{ for }(0,t)\in\Lambda.
\]
Note that $\Phi=\varphi$ inside the rectangle $R(\delta_1,\varepsilon_0)$, because $\vartheta=1$ for $|t|<\varepsilon_0$ (identifying $t\in(-\varepsilon_0,\varepsilon_0)$ with its class). Thus, to prove our statement, we only have to lift this function $\tilde{\pi}$ from the neighbourhood $W$ of the circle $\mathbb{R}/\ell \mathbb{Z}$ to a neighbourhood of the line $\{z=0\}$ of the plane.

To this end, choose $c_1\le\min(\varepsilon_0,\delta_1)$ such that $B:=R(c_1,c_1)\subset W$. $B$ is simply connected and thus $\tilde\pi|_B$ lifts to a smooth $\pi:B\to\R$ with $\pi(0,0)=0$. The function $(\pi(0,t)-t)/\ell$ is continuous, having integer values and vanishes at $t=0$, so $\pi(0,t)=t$. Now choose a $c\le c_1$ such that, for the rectangle $B':=R(c,c)$, it is $\varphi(B')\subset B$. On $B'$, $\Phi=\varphi$, so $\pi\circ\varphi\equiv\pi\pmod\ell$; hence $(\pi\circ\varphi-\pi)/\ell$ is continuous, with integer values on the connected set $B'$ and vanishes at $0$; therefore $\pi\circ\varphi=\pi$ on $B'$.

\end{proof}

The function $\pi$, given in the last proposition, is the sought-after function $\pi_2$ and we now proceed to define the function $\pi_1$.

\subsection{\texorpdfstring{The function $\pi_1$}{The function pi1}}

By \ref{prop:A}, there exists a smooth function $\pi$ near $0$ with $\pi(0,t)=t$. We use it to introduce new coordinates, in a neighbourhood of the origin.

Note that the map $\chi(z,t):=(z,\pi(z,t))=:(z,p)$ has Jacobian determinant $\pi_t$ and $\pi_t(0,0)=1$, which means that $\chi$ is a diffeomorphism from a neighbourhood of $0$ onto a rectangle
$Q=(-\delta,\delta)\times(-\eta,\eta)$. This diffeomorphism preserves $Z=\{z=0\}$ and is the identity on it, since $\chi(0,t)=(0,t)$. Also, $Q$ is a rectangle and thus $(z,p)\in Q$ implies $(sz,p)\in Q$ for $s\in[0,1]$.

From the relation $\dd z\wedge\dd\pi=\pi_t\,\dd z\wedge\dd t$ we get (for points not in $Z$):
\[
\omega=\frac{\dd z\wedge\dd t}{z}=\sigma\,\frac{\dd z\wedge\dd\pi}{z},\quad\text{where } \sigma:=\frac1{\pi_t},\quad \sigma(0,t)=1 .
\]
Define the function $\hat\sigma:=\sigma\circ\chi^{-1}$ and, in the coordinates $(z,p)$, we have:
\[
\omega=\hat\sigma(z,p)\,\frac{\dd z\wedge\dd p}{z},\qquad \hat\sigma(0,p)=1 .
\]
By Hadamard's lemma, we can write $\hat\sigma=1+z\,\hat r$ where the function $\hat r(z,p)=\int_0^1\partial_z\hat\sigma(sz,p)\,\dd s$ is smooth on $Q$ and hence:
\begin{equation}\label{eq:omega-zp}
\omega=\frac{\dd z\wedge\dd p}{z}+\hat r\,\dd z\wedge\dd p .
\end{equation}

We are now looking for a function $\lambda(z,p)$, such that the mapping $
\Theta(z,p):=(z\,e^{\lambda(z,p)},\,p)$ pulls the form $\omega_{s}=\frac{\dd w}{w}\wedge\dd p$ back to $\omega$. Since $\dd w/w=\dd z/z+\dd\lambda$, we compute that:
\[
\Theta^*\omega_{s}=\Bigl(\frac{\dd z}{z}+\dd\lambda\Bigr)\wedge\dd p=\frac{\dd z\wedge\dd p}{z}+\lambda_z\,\dd z\wedge\dd p .
\]
Comparing this with equation \eqref{eq:omega-zp}, we see that the condition $\Theta^*\omega_{s}=\omega$ holds if:
\[
\partial_z\lambda=\hat r,\quad\text{which has the solution } \lambda(z,p):=\int_0^z\hat r(s,p)\,\dd s .
\]

Composing with $\chi$, we get the smooth map:
\[
\Psi(z,t)=\Theta \circ \chi\,(z,t)=\Bigl(z\,\pi_1(z,t),\ \pi_2(z,t)\Bigr),\quad \pi_1(z,t):=\exp\!\int_0^z\hat r\bigl(s,\pi(z,t)\bigr)\,\dd s .
\]

\begin{lemma}\label{lem:a}
Let $\Psi:=(w,\pi_2)$, where $w:=z\,\pi_1$, be the map constructed above. Then, in a neighbourhood of $0$, $\Psi\in\Gb$ and $\rho(\Psi)=1$.
\end{lemma}

\begin{proof}
Since $\pi_1=e^{\lambda\circ\chi}>0$, we have $w=0$ if and only if $z=0$. For $z=0$, $\pi(0,t)=t$ and $\lambda(0,\cdot)=0$ give $\Psi(0,t)=(0,t)$. At $(0,t)$:
\[
w_z=\pi_1+z\,\partial_z\pi_1=\pi_1(0,t)=1,\qquad w_t=z\,\partial_t\pi_1=0,\qquad \pi_t=1.
\]
Thus, the Jacobian matrix at the origin has determinant equal to $1$.

It is $\Psi=\Theta\circ\chi$ with $\Theta(z,p)=(ze^{\lambda},p)$, and $\Theta$ pulls $\frac{\dd w}{w}\wedge\dd p$ back to:
\[
\Bigl(\frac{\dd z}{z}+\dd\lambda\Bigr)\wedge\dd p=\Bigl(\frac1z+\hat r\Bigr)\dd z\wedge\dd p=\hat\sigma\,\frac{\dd z\wedge\dd p}{z}.
\]
Pulling back by $\chi$ ($z\mapsto z$, $p\mapsto\pi$):
\[
\Psi^*\Bigl(\frac{\dd w}{w}\wedge\dd\pi\Bigr)=\sigma\,\frac{\dd z\wedge\dd\pi}{z}=\frac{\dd z\wedge\dd t}{z},
\]
the last equality being the identity $\omega=\sigma\frac{\dd z\wedge\dd\pi}{z}$ established above.

Also $\rho(\Psi)=\pi_1(0,\cdot)=e^{\lambda(0,\cdot)}=1$ and, therefore, $\Psi \in \Gb_1$.
\end{proof}

\subsection{Proof of the theorem}

We are now ready to prove the main theorem of this section.

\begin{proof}[Proof of \ref{thm:A}]

Let $\Psi=(w,\pi_2)$, $w=z\pi_1$, as given by \ref{lem:a} and define:
\[
\varphi':=\Psi\circ\varphi\circ\Psi^{-1}\in\Gb.
\]
Note that $\rho(\varphi')=\rho(\varphi)=u_0$.

Let $y=(w,\tau)$ be a point near the origin and put $x:=\Psi^{-1}(y)$. Then $\Psi(x)=y$ and, since $\Psi=(w,\pi)$, we get $w(x)=w$ and $\pi(x)=\tau$.

We also have:
\[
\varphi'(y)=\Psi\bigl(\varphi(\Psi^{-1}(y))\bigr)=\Psi\bigl(\varphi(x)\bigr)=\bigl(w(\varphi(x)),\ \pi(\varphi(x))\bigr),
\]
and, since $\pi\circ\varphi=\pi$, $\pi(\varphi(x))=\pi(x)=\tau$.
Hence:
\[
\varphi'(w,\tau)=\bigl(F'(w,\tau),\ \tau\bigr),\qquad \text{where } F'(w,\tau):=w(\varphi(\Psi^{-1}(w,\tau))).
\]
That is, the second component of $\varphi'$ is $T'(w,\tau)=\tau$. Since $\varphi'$ is also a member of $\Gb$, we can write $F'=w\,u'$, with $u'$ smooth. The Jacobian equation, for $T'=\tau$, reads as $F'_w=u'$, i.e.\ $u'+w\,u'_w=u'$. Thus $w\,u'_w=0$, so $u'_w=0$ off $Z$, hence everywhere by continuity. Therefore $u'(w,\tau)=u'(0,\tau)=u_0(\tau)$ (on a rectangle, integrating in $w$), and:
\[
\varphi'=L_{u_0}.
\]
We have proved that $\Psi\varphi\Psi^{-1}=L_{u_0}$, with $\Psi\in\Gb_1$.
\end{proof}
Thus, every conjugacy class of normally hyperbolic b-symplectomorphisms is completely determined by a member of $\Ab$. To be more precise, we have the following:

\begin{corollary}\label{cor:invariant}
There exists a bijection between conjugacy classes of normally hyperbolic b-symplectomorphisms and functions $a\in\Ab$ with $|a(0)|\neq1$. The normally hyperbolic germs $\varphi,\varphi'\in\Gb$ belong to the same conjugacy class if, and only if, $\rho(\varphi)=\rho(\varphi')$.
\end{corollary}

We now turn our attention to the case of vector fields preserving the form $\omega$.

\section{Local study of b-symplectic vector fields}\label{sec:vf}

\subsection{b-symplectic vector fields}

A $b$-vector field on an open set $\Omega\subset\R^2$ is a smooth vector field tangent to $Z$ and it is thus necessarily of the form $v=A\,z\partial_z+B\,\partial_t$, where $A,B\in C^\infty(\Omega)$. On $\Omega\setminus Z$ it is:
\begin{equation}\label{eq:iota}
\iota_v\omega=A\,\dd t-B\,\frac{\dd z}{z}.
\end{equation}
We shall call the vector field $v$ ``\emph{$b$-symplectic}" if $\Lie_v\omega=0$ on $\Omega\setminus Z$. Since $\dd\omega=0$ (it is a $3$-form on a surface), Cartan's magic formula gives us $\Lie_v\omega=\dd(\iota_v\omega)$, which means that $v$ is $b$-symplectic if, and only if, $\iota_v\omega$ is closed.

On the other hand, for $f\in C^\infty(\Omega)$, define the vector field:
\[
X_f:=z\,f_t\,\partial_z-z\,f_z\,\partial_t .
\]
This is a smooth $b$-vector field on $\Omega$, which vanishes on $Z\cap\Omega$, with $\iota_{X_f}\omega=\dd f$ on $\Omega\setminus Z$ and it is the only smooth vector field with this property. Furthermore, $X_f=X_g$ if, and only if, $f-g$ is locally constant, and we note that not every $b$-symplectic field is of the form $X_f$. Indeed, the field $\partial_t$, for example, has $\iota_{\partial_t}\omega=-\frac{\dd z}{z}$, which is not the differential of a smooth function.

Non-vanishing b-symplectic vector can be straightened, via a b-symple-ctomorphism. 

\begin{theorem}\label{thm:flowbox}
Let $v$ be the germ at $0$ of a $b$-symplectic vector field with $v(0)\neq0$. Then a $\Psi\in\Gb$ and a constant $c\neq0$ exist such that $\Psi_*v=c\,\partial_t$.
\end{theorem}
\begin{proof}
We write $v=A\,z\partial_z+B\,\partial_t$ and thus $v(0)=B(0,0)\,\partial_t$. If $c:=B(0,0)$, then $c\neq0$ because $v(0)\neq0$. By Hadamard's lemma, the function $z\mapsto B(z,0)$ can be written as:
\[
B(z,0)=c+z\,B_1(z),\qquad \text{where }B_1\ \text{smooth}.
\]
Let $\varphi_s$ be the flow of $v$. Since $v$ is tangent to $Z$, this flow preserves $Z$, so $\varphi_s$ maps $\{z\neq0\}$ into itself. Since $\iota_v\omega$ is closed off $Z$, $\Lie_v\omega=0$ there. Hence $\varphi_s^*\omega=\omega$ off $Z$.

Now, define the map $\Theta_0(z,t):=\varphi_{t/c}(z,0)$ and note that it is a diffeomorphism, with $(\Theta_0)_*(c\,\partial_t)=v$.

If $z\neq0$ and $s=t/c$, we calculate that $\partial_z\Theta_0(z,t)=d_{(z,0)}\varphi_s\cdot \partial_z$ and $\partial_t\Theta_0(z,t)=\tfrac1c v(\varphi_s(z,0))$. Since $\varphi_s^*\omega=\omega$, we get:
\[
\Theta_0^*\omega=\frac{B(z,0)}{c\,z}\,\dd z\wedge\dd t,
\]
and thus:
\[
\Theta_0^*\omega=\Bigl(\frac1z+\beta(z)\Bigr)\dd z\wedge\dd t,\qquad \beta:=\frac{B_1}{c}.
\]

Let $\zeta(z):=z\exp\int_0^z\beta$. This is a diffeomorphism of $(\R,0)$, with $\zeta'(0)=1$, and:
\[
\frac{\zeta'}{\zeta}=\frac1z+\beta .
\]
Now, define $\Theta(z,t):=\Theta_0\bigl(\zeta^{-1}(z),t\bigr)$. Then $\Theta$ is a diffeomorphism germ with $\Theta_*(c\,\partial_t)=v$ (only the $z$-variable was reparametrised). For points not belonging in $Z$, we have:
\[
\Theta^*\omega=\frac{\zeta'}{\zeta}\bigl(\zeta^{-1}(z)\bigr)\,(\zeta^{-1})'(z)\,\dd z\wedge\dd t=\frac{1}{\zeta(\zeta^{-1}(z))}\,\dd z\wedge\dd t=\omega .
\]
We thus have $\Theta\in\Gb$ and the map $\Psi:=\Theta^{-1}\in\Gb$ satisfies $\Psi_*v=c\,\partial_t$.
\end{proof}

Now, we consider b-symplectic vector fields that vanish at the origin.

\begin{lemma}\label{lem:hamiltonian}
Let $v=A\,z\partial_z+B\,\partial_t$ be a $b$-symplectic vector field on the rectangle $R=R(\delta,\varepsilon)$, with $v(0)=0$. Then $B(0,t)=0$ for all $|t|<\varepsilon$ and there exists a smooth function $f:R\to\R$, unique up to an additive constant, such that $v=X_f$. In particular $A(0,t)=f_t(0,t)$.
\end{lemma}
\begin{proof}
Consider the 1-form  $\eta=A\,\dd t-\frac Bz\,\dd z$, defined on $R\setminus Z$. It is closed, since $v$ is $b$-symplectic. We compute that $\dd\eta=\bigl(A_z+\frac{B_t}{z}\bigr)\dd z\wedge\dd t$, so we must have $zA_z+B_t=0$ on $R\setminus Z$, hence everywhere in $R$ by continuity. For $z=0$, this equation gives us $B_t(0,t)=0$ and, as $v(0)=(0,B(0,0))=0$, we get $B(0,t)=0$ for all $t$. Using Hadamard's lemma, we write $B=zB_1$ with $B_1$ smooth, and $\eta=A\,\dd t-B_1\,\dd z$ extends to a smooth, closed, $1$-form everywhere in the rectangle $R$; we can, therefore, write $\eta=\dd f$ by the Poincar\'e lemma: $A=f_t$, $B_1=-f_z$. We get $v=X_f$, while uniqueness up to constants is clear.
\end{proof}

The flow of such a vector field consists of b-symplectomorphisms.

\begin{lemma}\label{lem:flow}
Let $v=X_f$ on $R=R(\delta,\varepsilon)$ and $\lambda(t):=f_t(0,t)$. There is a neighbourhood $U_1$ of the origin such that the flow $\varphi_s$ of $v$ is defined on $U_1$ for $|s|\le1$ and:
\begin{enumerate}[label=(\roman*)]
\item $f\circ\varphi_s=f$;
\item $\varphi_s\in\Gb$ for every $s$, and $\rho(\varphi_s)=e^{s\lambda}$.
\end{enumerate}
\end{lemma}
\begin{proof}
Existence of $U_1$ follows from $v(0)=0$ and continuous dependence on initial conditions.

(i) $\dd f(v)=(\iota_v\omega)(v)=\omega(v,v)=0$ on $R\setminus Z$, hence everywhere on $R$ by continuity.

(iii) Points of $Z\cap R$ are zeros of $v$, hence stationary, and a phase curve with initial condition not on $Z$ never meets $Z$. We get that $O_{\varphi_s}=U_1\setminus Z$ and $\frac{\dd}{\dd s}\varphi_s^*\omega=\varphi_s^*\Lie_v\omega=0$, so $\varphi_s^*\omega=\omega$ there. Obviously $\varphi_s$ is a diffeomorphism onto its image, having as inverse the mapping $\varphi_{-s}$, while $\varphi_s(0)=0$, thus $\varphi_s \in \Gb$.

For the last assertion, let us write $\varphi_s=(F_s,T_s)$. The $z$-component of the flow equation is:
\[
\partial_sF_s=A(F_s,T_s)\,F_s ,\qquad F_0(z,t)=z .
\]
Differentiating this relation with respect to $z$, we get:
\[
\partial_s\partial_zF_s=\partial_zF_s\cdot A(F_s,T_s)+F_s\bigl[A_z(F_s,T_s)\partial_zF_s+A_t(F_s,T_s)\partial_zT_s\bigr].
\]
At the point $(0,t)\in U_1$ we have $F_s=0$ and $T_s=t$ and thus the bracket is multiplied by $0$. The function $y(s):=\partial_zF_s(0,t)$ satisfies the equation $y'=\lambda(t)y$, with $y(0)=1$. Hence $y(s)=e^{s\lambda(t)}$, which gives us $u_0$ by \ref{lem:structure}.
\end{proof}

The group $\Gb$ acts on the space of b-symplectic vector fields $X_g$ of interest by the rule $(\Psi,X_g)\mapsto \Psi\cdot X_g:=\Psi_*X_g$ and we wish to study the orbits of this action. As it turns out, the functions of all the b-symplectic vector fields belonging to the same orbit are related to each other.

\begin{lemma}\label{lem:equiv}
If $\Psi\in\Gb$, $\Psi:U\to V$ and $g\in C^\infty(U)$, then $\Psi_*X_g=X_{g\circ\Psi^{-1}}$ on $V$.
\end{lemma}
\begin{proof}
Just compute that:
\[
\iota_{\Psi_*X_g}\omega=(\Psi^{-1})^*\bigl(\iota_{X_g}\Psi^*\omega\bigr)=(\Psi^{-1})^*\dd g=\dd(g\circ\Psi^{-1})=\iota_{X_{g\circ\Psi^{-1}}}\omega.
\]
\end{proof}

Since the classification of the vector fields of interest reduces to a classification of functions, we turn our attention to functions defined on the b-symplectic plane.

\subsection{A classification of functions at regular points}

Let $\mathfrak m$ be the ring of germs at the origin of smooth functions vanishing at $0\in \mathbb{R}^2$. \ref{lem:hamiltonian} ensures us that $f\mapsto X_f$ is a bijection from $\mathfrak m$ onto the $b$-symplectic fields vanishing at $0$, and by \ref{lem:equiv}, the action $\Psi\cdot f=f\circ\Psi^{-1}$ of $\Gb$ on $\mathfrak{m}$ gives us the classification of $b$-symplectic fields, vanishing at $0$, up to members of $\Gb$.

Some first invariants are easily obtained.
 
\begin{proposition}\label{prop:f0}
Let $f,g\in\mathfrak m$ and suppose $g=f\circ\Psi^{-1}$ for some $\Psi\in\Gb$, that is, $X_f$ and $X_g$ are $\Gb$-conjugate. Then $g|_Z=f|_Z$ as germs at $0\in\R$. In particular, $f_0(t):=f(0,t)$ and $\lambda(t)=f_0'(t)$ are invariants of the conjugacy class of $X_f$.
\end{proposition}
\begin{proof}
$\Psi(0,t)=(0,t)$ for all small $|t|$ and thus $g(0,t)=f\bigl(\Psi^{-1}(0,t)\bigr)=f(0,t)$. Hence $f_0=g_0$ and therefore $\lambda_f=f_0'=g_0'=\lambda_g$.
\end{proof}

To proceed with our study, we shall denote by $\Dd _f$ the differential operator $f_z\,\partial_t-f_t\,\partial_z$. In this way, $X_f=-z\,\Dd_f$ and $\forall f,h\in C^\infty$, it is $X_h(f)=z\,\Dd_f(h)$. Also, $\Omega^k$ stands for the space of germs of smooth $k$-forms, so $\Omega^0=C^\infty$ and $\Omega^2=C^\infty\,\dd z\wedge\dd t$, while $\langle z\rangle=\{zg:g\in C^\infty\}$ the ideal generated by $z$ and $\mathfrak m_t$ is the space of germs at $0\in\R$ of smooth functions of $t$ vanishing at $0$.

A smooth curve, passing through the identity, in $\Gb$ is a family $\{\Psi_s\}_{|s|<s_0}$ of members of $\Gb$, for some $s_0>0$, with $\Psi_0=\mathrm{Id}$. For $f\in\mathfrak m$, the tangent space of the $\Gb$-orbit $\mathcal{O}_f$ of $f$, at the point $f$, is:
\[
T_f\mathcal{O}_f:=\Bigl\{\tfrac{\dd}{\dd s}\Big|_{s=0}\,f\circ\Psi_s^{-1}\ :\ \{\Psi_s\}\text{ is a smooth curve in }\Gb\text{ with }\Psi_0=\mathrm{Id}\Bigr\}\subset C^\infty .
\]

\begin{proposition}\label{prop:Tf}
Let $f\in\mathfrak m$.
\begin{enumerate}[label=(\roman*)]

\item $T_f\mathcal{O}_f=\{X_h(f)/h\in C^\infty\}=z\cdot\Dd_f(C^\infty)\subset\langle z\rangle$.

\item Multiplication by $z$ and by $\dd z\wedge\dd t$ induce linear isomorphisms:
\[
C^\infty/\Dd_f(C^\infty)\ \xrightarrow{\ \cong\ }\ N(f):=\langle z\rangle/T_f\mathcal{O}_f
\]
\[
C^\infty/\Dd_f(C^\infty)\ \xrightarrow{\ \cong\ }\ \Omega^2/\bigl(\dd f\wedge\dd\Omega^0\bigr),
\]
where $\dd f\wedge\dd\Omega^0:=\{\dd f\wedge\dd h:h\in C^\infty\}$.
\item Let $r:\mathfrak m\to\mathfrak m_t$, $r(g)(t)=g(0,t)$. Then $r$ is surjective, $\ker r=\langle z\rangle$, and the sequence:
\[
0\longrightarrow N(f) \xrightarrow{\ \iota\ }\mathfrak m/T_f\mathcal{O}_f\xrightarrow{\ \bar r\ }\mathfrak m_t\longrightarrow0
\]
is exact, with $\iota$ induced by the inclusion $\langle z\rangle \subset\mathfrak m$ and $\bar r$ induced by $r$.
\item $T_f\mathcal{O}_f$ has infinite codimension in $\mathfrak m$, and the $\Gb$-orbit of $f$ lies in $f+\langle z\rangle$.
\end{enumerate}
\end{proposition}

\begin{proof}
\textbf{(i)} Define the vector field $V:=\partial_s \Psi_s|_{0}$. We shall show that it is b-symplectic.

Since $\Psi_s(0,t)=(0,t)$, we get $V(0,t)=0$ by differentiating with respect to $s$. Thus $V$ is a smooth field which vanishes on $Z$ and, using  Hadamard's lemma, $V=zA\,\partial_z+B\,\partial_t$ with $A,B$ smooth and $B(0,t)=0$, i.e.\ $V$ is a $b$-vector field with $V(0)=0$.

Note also that $\Psi_s^*\omega$ depends smoothly on $(s,x)$, and $\partial_s\Psi_s^*\omega|_{0}=\Lie_V\omega$. Due to the fact that $\Psi_s^*\omega=\omega$ for all $s$, the left-hand side of the last equation  is $0$ and thus we get $\Lie_V\omega=0$, off $Z$.

We have proved that $V$ is a b-symplectic vector field and thus \ref{lem:hamiltonian} ensures us that $V=X_h$, for some $h\in C^\infty(R)$. The velocity of the curve $\{\Psi_s\}$ at the point $f$ is therefore $-X_h(f)=-z\,\Dd_f(h)=z\,\Dd_f(-h)$. Hence every element of $T_f\mathcal{O}_f$ is of the form $X_{h'}(f)$, for $h'=-h$, and thus lies in $z\,\Dd_f(C^\infty)\subset \langle z\rangle$.

To prove the opposite inclusion, if $h\in C^\infty$, let $\Psi_s$ be the flow of $X_{h}$. We know that $\Psi_s\in\Gb$ and therefore:
\[
\frac{\dd}{\dd s}f\circ\Psi_s^{-1}\Big|_{0}=-X_{h}(f)=X_{-h}(f).
\]
Since $h$ ranges over a vector space, $\{-X_h(f)\}=\{X_h(f)\}$, which gives $\{X_h(f):h\in C^\infty\}\subset T_f\mathcal{O}_f$ and hence $T_f\mathcal{O}_f=z\,\Dd_f(C^\infty)$.

\medskip
\textbf{(ii)} Consider the surjective map $m_z:C^\infty\to \langle z\rangle$, $g\mapsto zg$. It is also injective: $zg=0$ implies $g=0$ on $\{z\ne0\}$, hence $g\equiv0$ by continuity. Since we have proved that $m_z(\Dd_f(C^\infty))=T_f\mathcal{O}_f$, $m_z$ induces an isomorphism $C^\infty/\Dd_f(C^\infty)\cong(z)/T_f\mathcal{O}_f$.

Now, consider the linear isomorphism $\mu:C^\infty\to\Omega^2$, $g\mapsto g\,\dd z\wedge\dd t$. For $f,h\in C^\infty$:
\[
\dd f\wedge\dd h=(f_zh_t-f_th_z)\,\dd z\wedge\dd t=\Dd_f(h)\,\dd z\wedge\dd t .
\]
Thus $\mu(\Dd_f(h))=\dd f\wedge\dd h$ and so $\mu(\Dd_f(C^\infty))=\dd f\wedge\dd\Omega^0$; $\mu$ induces therefore an isomorphism $C^\infty/\Dd_f(C^\infty)\cong\Omega^2/(\dd f\wedge\dd\Omega^0)$.

\medskip
\textbf{(iii)} If $g\in\mathfrak m$, we have $r(g)(0)=g(0,0)=0$ and so $r(g)\in\mathfrak m_t$. Given $\phi\in\mathfrak m_t$, the function $g(z,t)=\phi(t)$ lies in $\mathfrak m$ and satisfies $r(g)=\phi$. The map $r$ is well defined and onto.

We now compute the kernel of $r$. If $g\in\mathfrak m$ and $g(0,t)=0$, then by Hadamard's lemma:
\[
g(z,t)=\int_0^1\frac{\dd}{\dd\sigma}g(\sigma z,t)\,\dd\sigma=z\int_0^1g_z(\sigma z,t)\,\dd\sigma\in \langle z\rangle.
\]
Conversely, every element of $\langle z\rangle$ vanishes on $Z$ and thus lies both in $\mathfrak m$ and in $\ker r$. Hence $\ker r=\langle z\rangle$.

Since $T_f\mathcal{O}_f\subset \langle z\rangle \subset\mathfrak m$, the map:
\[\iota:\langle z\rangle /T_f\mathcal{O}_f\to\mathfrak m/T_f\mathcal{O}_f,\quad g+T_f\mathcal{O}_f\mapsto g+T_f\mathcal{O}_f,
\]
is well defined, and it is injective, because $g\in T_f\mathcal{O}_f$ whenever $g+T_f\mathcal{O}_f=0$. The map $\bar r(g+T_f\mathcal{O}_f)=r(g)$ is well defined (since $r(T_f\mathcal{O}_f)=0$), and is surjective since $r$ is. Finally:
\[
\ker\bar r=\{g+T_f\mathcal{O}_f: r(g)=0\}=\langle z\rangle/T_f\mathcal{O}_f=\operatorname{im}\iota ,
\]
which proves the exactness of the sequence.

\medskip
\textbf{(iv)} By (iii), $\mathfrak m/T_f\mathcal{O}_f$ surjects onto $\mathfrak m_t$, which is infinite dimensional and hence $\operatorname{codim}_{\mathfrak m}T_f\mathcal{O}_f=\infty$. If $g=f\circ\Psi^{-1}$ with $\Psi\in\Gb$, then $g|_Z=f|_Z$ \ref{prop:f0}, i.e.\ $r(g-f)=0$, so $g-f\in \langle z\rangle$. Thus the orbit lies in $f+\langle z\rangle$, and $\operatorname{codim}_{\mathfrak m}\langle z\rangle=\dim\mathfrak m_t=\infty$.
\end{proof}

Since $T_f\mathcal{O}_f$ is always of infinite codimension in $\mathfrak{m}$, we study the space $N_f$ instead, the structure of which depends on whether the origin is a critical point for $f$ or not.

\begin{theorem}\label{thm:dichotomy}
Let $f\in\mathfrak m$. If $\dd_0 f\neq0$, it is $N(f)=0$, while if $\dd_0 f=0$, $\dim N(f)=\infty$.
\end{theorem}
\begin{proof}
Let us suppose that $\dd_0 f\neq0$. The vector field $\Dd_f=f_z\,\partial_t-f_t\,\partial_z$ does not vanish at the origin and, by the flow-box theorem, there exists a diffeomorphism $\Phi=(\xi,\tau)$ from a neighbourhood of $0$ onto a box $B=I\times J$ (with $J\ni0$ an interval), such that $\Phi_*\Dd_f=\partial_\tau$.

For a smooth germ $q$, define $\tilde q=q\circ\Phi^{-1}$ and:
\[
\tilde h(\xi,\tau)=\int_0^\tau\tilde q(\xi,\sigma)\,\dd\sigma .
\]
In this way, $h:=\tilde h\circ\Phi$ is a smooth germ with:
\[
\Dd_f h=(\partial_\tau\tilde h)\circ\Phi=\tilde q\circ\Phi=q,
\]
therefore $\Dd_f$ is surjective. By \ref{prop:Tf}(i):
\[
T_f\mathcal{O}_f=z\,\Dd_f(C^\infty)=zC^\infty=\langle z\rangle,
\]
and therefore $N(f)=\langle z\rangle/T_f\mathcal{O}_f=0$.

Let us now suppose that $\dd_0 f=0$, which means that $f_z,f_t\in\mathfrak m$. If  $g\in\mathfrak m^{k+1}$, we have that $\Dd_fg=f_z\,g_t-f_t\,g_z\in\mathfrak m^{k+1}$.

Hence $\Dd_f(\mathfrak m^{k+1})\subset\mathfrak m^{k+1}$, and thus $\Dd_f$ induces a linear map:
\[
\bar\Dd_k:J^k\to J^k,\qquad  \bar\Dd_k(g+\mathfrak m^{k+1})=\Dd_fg+\mathfrak m^{k+1},
\]
where $J^k$ is the space of $k$-jets $J^k:=C^\infty/\mathfrak m^{k+1}$, with $\dim J^k=\binom{k+2}{2}<\infty$.
 
Consider now the composition:
\[
C^\infty\xrightarrow{\ \pi\ }J^k\longrightarrow J^k/\bar\Dd_k(J^k),
\]
where $\pi$ is the quotient map. It is surjective and $\Dd_f(C^\infty)$ belongs to its kernel, because $\pi(\Dd_fg)=\bar\Dd_k(\pi g)\in\bar\Dd_k(J^k)$.
Therefore it induces a surjection:
\[
N(f)\cong C^\infty/\Dd_f(C^\infty)\twoheadrightarrow J^k/\bar\Dd_k(J^k),
\]
and so $\dim N(f)\ge\dim\bigl(J^k/\bar\Dd_k(J^k)\bigr)$, for every $k$. Hence, we have that:
\[
\dim N(f)\ \ge\ \dim\ker\bar\Dd_k,\quad\ \text{for all } k
\]
and it remains to show that $\dim\ker\bar\Dd_k\to\infty$, as $k\to\infty$.

Since $\Dd_f(f)=0$, the chain rule gives us $\Dd_f(f^j)=jf^{j-1}\Dd_f(f)=0$, for every $j\ge0$. So $f^j\in\ker\Dd_f$, and therefore its class $[f^j]_k$ lies in $\ker\bar\Dd_k$.

We distinguish two cases, according to the Taylor series of $f$ at $0$.

If the Taylor series of $f$ is not zero, let $r$ be the order of $f$, that is, the lowest degree of a nonzero homogeneous term of
its Taylor series. Then $r\ge2$ because $f(0)=0$ and $\dd_0f=0$. Write $f_r$ for that homogeneous polynomial of degree $r$. The lowest-degree
term of $f^j$ is $f_r^{\,j}$, which is nonzero because the polynomial ring has no zero divisors. So $f^j$ has order exactly $rj$.
Consider $0\le j\le\lfloor k/r\rfloor$, so that $rj\le k$. Then $[f^j]_k\neq0$ in $J^k$. These classes are linearly independent: if
$\sum_jc_j[f^j]_k=0$ with not all $c_j=0$, let $j_0$ be the smallest index with $c_{j_0}\neq0$. Then $\sum_jc_jf^j$ has order exactly
$rj_0\le k$ (the lowest-degree term is $c_{j_0}f_r^{\,j_0}$, and the others have higher degree), so it is not in $\mathfrak m^{k+1}$, a contradiction.
Therefore:
\[
\dim\ker\bar\Dd_k\ \ge\ \lfloor k/r\rfloor+1\ \xrightarrow[k\to\infty]{}\ \infty .
\]

If the Taylor series of $f$ is zero (i.e.\ $f$ is flat), then $f_z,f_t$ are also flat, so for every $g$ we have:
$\Dd_fg=f_zg_t-f_tg_z\in\bigcap_k\mathfrak m^{k+1}$. Thus $\bar\Dd_k=0$ and $\ker\bar\Dd_k=J^k$, with
$\dim J^k=\binom{k+2}{2}\to\infty$.

In both cases, since $\dim N(f)\ge\dim\ker\bar\Dd_k$ for all $k$ and the right-hand side is unbounded, we get $\dim N(f)=\infty$.
\end{proof}

To state the main classification result, we need two lemmata. 

\begin{lemma}\label{lem:param}
Let $J$ be an open interval containing $[0,1]$ and $W_0$ a neighbourhood of $0\in \mathbb{R}^2$. Let also $g_s:W_0\to \mathbb{R}$ be a family of smooth germs, depending smoothly on $s\in J$ and $\Sigma$ a straight line through $0$ with $\Dd_{g_s}(0)\notin T\Sigma$, for all $s\in[0,1]$. Then there exists a neighbourhood $W$ of $0$ such that, for every smooth family of functions $q_s$, there is another smooth family $h_s$, defined on $W$, with $\Dd_{g_s}h_s=q_s$, $s\in[0,1]$.
\end{lemma}
\begin{proof}
Let $e$ span $T\Sigma$ and $\Phi^s_\tau$ be the flow of $V_s=\Dd_{g_s}$. The map $(s,\xi,\tau)\mapsto(s,\Phi^s_\tau(\xi e))$ is smooth, and its differential at $(s,0,0)$ is invertible because $V_s(0)\notin\R e$. Since it is a local diffeomorphism at each point of the compact set $[0,1]\times\{(0,0)\}$, on which it is injective, it is injective on a neighbourhood of that set, hence a diffeomorphism onto an open set containing $[0,1]\times W$ for some neighbourhood $W$ of $0$. In these coordinates, $V_s=\partial_\tau$ and $h_s(\Phi^s_\tau(\xi e)):=\int_0^\tau q_s(\Phi^s_{\tau'}(\xi e))\,\dd\tau'$ solves the equation.
\end{proof}

\begin{lemma}\label{lem:moser}
Let $g_s$, $s\in J\supset[0,1]$, be a smooth family on $J\times W_0$ satisfying the hypotheses of \ref{lem:param} (i.e.\ $\Dd_{g_s}(0)\notin T\Sigma$ for all $s\in[0,1]$). Suppose also that $\partial_sg_s=z\,q_s$, for a family $q_s$, smooth in $(s,z,t)$. Then there exist a neighbourhood $U$ of $0$ and a family $\Psi_s:U\to\R^2$, $s\in[0,1]$, smooth in $(s,x)$, satisfying $\Psi_0=\id$ and $\Psi_s\in\Gb$ for every $s$, such that:
\[
g_s\circ\Psi_s=g_0\quad\text{on }U,\qquad s\in[0,1].
\]
\end{lemma}

\begin{proof}
Applying \ref{lem:param} to the family $g_s$ and the family $-q_s$, we obtain a neighbourhood $W$ of $0$ and a smooth family $h_s$, defined on $[0,1]\times W$, such that:
\begin{equation}\label{eq:homological}
\Dd_{g_s}h_s=-q_s\quad\text{on }W,\ s\in[0,1].
\end{equation}
Let $X_s:=X_{h_s}$ and denote by $\Psi_s$ its flow. As we saw in \ref{lem:flow}, this flow consists of members of $\Gb$.

Fix $x$ and calculate that:
\begin{align*}
\frac{\dd}{\dd s}\bigl(g_s\circ\Psi_s\bigr)(x)
&=(\partial_sg_s)\bigl(\Psi_s(x)\bigr)+\dd g_s\bigl(\partial_s\Psi_s(x)\bigr)
=\bigl(\partial_sg_s+X_s(g_s)\bigr)\bigl(\Psi_s(x)\bigr)=\\
&= \Big( \partial_sg_s+X_{h_s}(g_s)\Big)\bigl(\Psi_s(x)\bigr)=\Big(z\,q_s+z\,\Dd_{g_s}(h_s)\Big)\bigl(\Psi_s(x)\bigr)=\\
&=\Big(z\,q_s-z\,q_s\Big)\bigl(\Psi_s(x)\bigr)=0 .
\end{align*}
Hence $\frac{\dd}{\dd s}(g_s\circ\Psi_s)(x)=0$ for all $s\in[0,1]$ and so $g_s\circ\Psi_s=g_0\circ\Psi_0=g_0$.
\end{proof}

In the space $\mathfrak{m}$ we define two strata:
\begin{enumerate}
\item The stratum $S_0$ consisting of those $f$ with $\partial_tf(0)\neq 0$.
\item The stratum $S_1$, consisting of those germs $f$ with  $\partial_tf(0)= 0$ and $\partial_zf(0)\neq 0$.
\end{enumerate}

We can now give models for the members of these two strata.

\begin{theorem}\label{thm:regular}
Let $f\in\mathfrak m$ and $f_0=f|_Z$.
\begin{enumerate}[label=(\alph*)]
\item If $f$ belongs to the stratum $S_0$, there exists a $\Psi\in\Gb_1$ with $f\circ\Psi^{-1}(z,t)=f_0(t)$. In this case, the $\Gb$-orbit of $f$ is $\{g\in\mathfrak m:\ g|_Z=f_0\}$, and $\Psi_*X_f=X_{f_0}=f_0'(t)\,z\partial_z$.
\item If $f$ belongs to the stratum $S_1$, there exists a $\Psi\in\Gb$ with $f\circ\Psi^{-1}(z,t)=f_0(t)+z$. In this case, the orbit of the germ $f$ is $\{g\in\mathfrak m:\ g|_Z=f_0,\ g_z(0)\neq0\}$, and $\Psi_*X_f=X_{f_0+z}=z\bigl(f_0'(t)\partial_z-\partial_t\bigr)$.
\end{enumerate}
\end{theorem}
\begin{proof}
Since $f-f_0(t)$ vanishes on $Z$, we can write:
\[
f(z,t)=f_0(t)+z\,r(z,t),\qquad r(z,t)=\int_0^1f_z(\sigma z,t)\,\dd\sigma ,
\]
with $r$ smooth. Then $r(0,0)=f_z(0)$, and $f_t(0)=f_0'(0)$.

(a) For the first assertion, put $g_s=f_0(t)+s\,z\,r$. Then the vector $\Dd_{g_s}(0)=-f_0'(0)\partial_z+s\,r(0)\partial_t$ is transverse to $\Sigma=Z$, for all $s$, and furthermore $\partial_sg_s=zr$.

\ref{lem:moser} provides us with a family $\Psi_s\in\Gb$, such that $g_s\circ\Psi_s=g_0=f_0$; take $\Psi=\Psi_1^{-1}$. The map $\Psi_1$ need not lie in $\Gb_1$ (recall that, as we saw in \ref{lem:flow}(iii), $\rho(\Psi_1)=\exp\int_0^1\partial_th_s(0,\cdot)\,\dd s$). To get $\Psi\in\Gb_1$ replace $\Psi$ by $L_{1/a}\circ\Psi$, $a=\rho(\Psi)$: since $f_0$ depends only on $t$, $f_0\circ L_{a}=f_0$, so the new map still satisfies $f\circ\Psi'^{-1}=f_0$, and $\rho(\Psi')=a^{-1}a=1$.

To compute the orbit of $f$, let $g\in\mathfrak m$ with $g|_Z=f_0$. Since $g_t(0)=f_0'(0)\neq0$, the germ $g$ is also in $S_0$ with the same restriction $f_0$ and the argument above shows that $g\sim f_0\sim f$. Conversely, if $g\sim f$, then $g|_Z=f_0$ by \ref{prop:f0}. Hence $\mathcal{O}_ f=\{g\in\mathfrak m: g|_Z=f_0\}$.

(b) For the constant $a=f_z(0)$, consider the map $L_a\in\Gb$ and define:
\[
\tilde f:=f\circ L_a^{-1}=f_0(t)+z\tilde r,
\]
where $\tilde r(0)=1$. If $g_s=f_0(t)+z+s\,z(\tilde r-1)$, then we have $\partial_zg_s(0)=1$ and $\partial_tg_s(0)=0$, which mean that the vector $\Dd_{g_s}(0)=\partial_t$ is always transverse to $\Sigma=\{t=0\}$ and $\partial_sg_s=z(\tilde r-1)$. \ref{lem:moser} then gives $\tilde f\sim f_0+z$, hence $f\sim f_0+z$. If $g|_Z=f_0$ and $g_z(0)\ne0$ the same argument applies, and $g_z(0)\neq0$ is preserved by $\Gb$ (it is multiplied by $1/u_0^\Psi(0)$ since $g_t(0)=0$).
\end{proof}

We see that we can classify the vector fields of interest, as long as they correspond to functions with a regular point at the origin. In this case, these functions provide us with the invariants of the conjugacy class. To be precise:

\begin{corollary}\label{cor:vfinv}
Let $f,g\in\mathfrak m$ with $\dd_0 f\neq0\neq \dd_0 g$. Then $X_f$ and $X_{g}$ are conjugate by an element of $\Gb$ if, and only if, $f_0'=g_0'$ as germs, where $f_0=f|_Z$ and $g_0=g|_Z$.
\end{corollary}
\begin{proof}
If $\Psi_*X_f=X_{g}$ then $g=f\circ\Psi^{-1}+c$, so $g_0=f_0+c$ (\ref{prop:f0}). Conversely, if $f_0'=g_0'$ then $f_0=g_0$ (since both of them vanish at $0$) and we get $f_t(0)=g_t(0)$. Both germs lie, thus, in the same stratum $S_0$ or $S_1$ (given $\dd f(0)\neq0\neq\dd g(0)$), and \ref{thm:regular} applies.
\end{proof}

\begin{remark} Here, we would like to point out the paper \cite{Kou}, in which the classification of functions defined on the Poisson plane $\{x,y\}=x$, under Poisson diffeomorphisms, is presented, in the complex holomorphic category (but see also \cite{Kou15} for pairs $(\omega,f)$ with $\omega$ degenerate along a curve). The Poisson structure considered there is the inverse of our $\omega$ (for $x=z$, $y=t$) and our \ref{thm:regular} corresponds to the normal forms presented in \cite{Kou} for the generic strata ($A_0$ and $A_1$ in Arnold's notation, the functional moduli being the germ $f_0$). The method used there was to fix $f$ in Arnold's normal form and classify the Poisson structures with respect to the diffeomorphisms preserving $f$ and $Z$; our approach is to fix $\omega$ and classify $f$, which can be done by Moser's method. We also note that our results hold in the real, $C^\infty$, category.
\end{remark}

\subsection{Orbits of 2-jets}

Having completed the classification in case the origin is a regular point for the germ $f$, we now turn our attention to the case where $f$ has a critical point at $0$. As we saw in \ref{thm:dichotomy}, in the second case the space complementary to $T_f\mathcal{O}_f$ is of infinite dimension and we expect more moduli than those appearing in the first case. This is already clear when studying the second differential of $f$.

Indeed, let $f\in\mathfrak m^2$. The $2$-jet of $f$ is the quadratic form:
\[
q(z,t)=a\,t^2+b\,zt+c\,z^2,\ \text{where} \quad a=\tfrac12f_{tt}(0),\quad b=f_{zt}(0),\quad c=\tfrac12f_{zz}(0).
\]
If $\Psi \in \Gb$, then the element $\Psi\cdot f=f\circ \Psi^{-1}$ depends only on the linear part of $\Psi$. By \ref{lem:structure}, it is $\Psi=(zu,\,t+zS)$ and so:
\[
D\Psi(0)=\begin{pmatrix}\alpha&0\\ \beta&1\end{pmatrix},\qquad \alpha=u_0(0)\neq0,\quad \beta=S(0),
\]
and $\Psi$ acts on $2$-jets only through the pair $(\alpha,\beta)$.

This leads to a classification of quadratic forms under a smaller group and  not under all of $\mathrm{GL}_2(\R)$. The substitutions allowed, as dictated by $D\Psi(0)$, are only the linear maps $z\mapsto z/\alpha$ and $t\mapsto t-\beta z/\alpha$, since they fix the line $Z=\{z=0\}$ pointwise. This is the reason that, in our classification, more invariants are found than the rank and the signature of ordinary quadratic forms. In particular $a=\tfrac12f_0''(0)$, the second derivative of $f|_Z$, is invariant (recall \ref{prop:f0}).

Regarding the orbits of 2-jets, we can prove the following:

\begin{proposition}\label{prop:jets}
Let $f\in\mathfrak m^2$, having as $2$-jet the quadratic form $q=a\,t^2+b\,zt+c\,z^2$. Put $\Delta:=b^2-4ac$.

\begin{enumerate}[label=(\alph*)]
\item For every $\alpha\neq0$ and $\beta\in\R$, there exists a $\Psi\in\Gb$, with $D\Psi(0)=\begin{pmatrix}\alpha&0\\ \beta&1\end{pmatrix}$.

\item If $\Psi\in\Gb$, with $D\Psi(0)=\begin{pmatrix}\alpha&0\\ \beta&1\end{pmatrix}$, then the $2$-jet of $f\circ\Psi^{-1}$ is $a\,t^2+b'\,zt+c'\,z^2$, with:
\[
b'=\frac{b-2a\beta}{\alpha},\qquad c'=\frac{a\beta^2-b\beta+c}{\alpha^2},\qquad \Delta'=b'^2-4ac'=\frac{\Delta}{\alpha^2}.
\]

\item $a$ and $\operatorname{sgn}\Delta$ are invariants of the $\Gb$-orbit of $f$. Moreover, if $a=0$, then $b=0$ is invariant, and if $a=b=0$, then $\operatorname{sgn}c$ is invariant.

\item The $2$-jet of $f$ is conjugate, by some $\Psi\in\Gb$, to exactly one of the models listed below.
\begin{center}
\begin{tabular}{@{}lll@{}}
\hline
type & conditions & model\\
\hline
(i) & $a\ne0$, $\Delta<0$ & $a\,t^2+\operatorname{sgn}(a)\,z^2$\\
(ii) & $a\ne0$, $\Delta>0$ & $a\,t^2-\operatorname{sgn}(a)\,z^2$\\
(iii) & $a\ne0$, $\Delta=0$ & $a\,t^2$\\
(iv) & $a=0$, $b\neq0$ & $zt$\\
(v) & $a=b=0$, $c\neq0$ & $\operatorname{sgn}(c)\,z^2$\\
(vi) & $a=b=c=0$ & $0$\\
\hline
\end{tabular}
\end{center}
\end{enumerate}
\end{proposition}

\begin{proof}
\begin{enumerate}[label=(\alph*)]
 
\item Just consider the map $(z,t)\mapsto (\alpha z,t+\beta z)$.

\item Write $A:=D\Psi(0)$. Since $\Psi^{-1}$ is smooth with $D\Psi^{-1}(0)=A^{-1}$ and $\Psi^{-1}(0)=0$, Taylor's formula gives $\Psi^{-1}(x)=A^{-1}x+O(|x|^2)$. Since $f=q+O(3)$ and $q$ is homogeneous of degree $2$, we get:
\[
f\bigl(\Psi^{-1}(x)\bigr)=q\bigl(A^{-1}x+O(|x|^2)\bigr)+O(|x|^3)=q(A^{-1}x)+O(|x|^3).
\]
So the $2$-jet of $f\circ\Psi^{-1}$ equals $q\circ A^{-1}$ and the rest are simple calculations.

\item By (b), the coefficient of $t^2$ is unchanged, so $a$ is invariant. Since $\Delta'=\Delta/\alpha^2$ and $\alpha^2>0$, $\operatorname{sgn}\Delta$ is invariant. If $a=0$, then $b'=b/\alpha$, so $b'=0$ if and only if $b=0$. If $a=b=0$, then $c'=c/\alpha^2$, so $\operatorname{sgn}c'=\operatorname{sgn}c$.

\item In every case, it suffices to choose $(\alpha,\beta)$ wisely.

\emph{Case $a\neq0$:} Take $\beta=b/(2a)$. Then $b'=0$, and:
\[
c'=\frac{1}{\alpha^2}\Bigl(\frac{b^2}{4a}-\frac{b^2}{2a}+c\Bigr)=-\frac{\Delta}{4a\,\alpha^2}.
\]
Hence $f\circ\Psi^{-1}$ has $2$-jet the form $a\,t^2+c'z^2$.
\begin{itemize}[nosep]
\item If $\Delta<0$, then $c'=\dfrac{|\Delta|}{4a\alpha^2}$ has the sign of $a$. Choosing $\alpha^2=\dfrac{|\Delta|}{4|a|}$ gives $c'=\operatorname{sgn}(a)$ and thus type (i).
\item If $\Delta>0$, then $c'=-\dfrac{\Delta}{4a\alpha^2}$ has the sign of $-a$. Choosing $\alpha^2=\dfrac{\Delta}{4|a|}$ gives $c'=-\operatorname{sgn}(a)$ and thus type (ii).
\item If $\Delta=0$, then $c'=0$ and we get type (iii).
\end{itemize}

\emph{Case $a=0$, $b\neq0$:} Then $b'=b/\alpha$ and $c'=(c-b\beta)/\alpha^2$. Take $\alpha=b$ and $\beta=c/b$. This gives us $b'=1$ and $c'=0$, i.e. type (iv).

\emph{Case $a=b=0$, $c\neq0$:} Take $\beta=0$. Then $b'=0$ and $c'=c/\alpha^2$. Choosing $\alpha=\sqrt{|c|}$ gives us $c'=\operatorname{sgn}(c)$ and type (v).

\emph{Case $q=0$:} Then $q\circ A^{-1}=0$ for every $A$ and we get type (vi).

To show that the six types are mutually inequivalent, we just use the invariants obtained in (c).
\end{enumerate}
\end{proof}

\begin{remark}
We have classified germs $f\in\mathfrak m$ under the action of $\Gb$. Arnold classified the same germs under the much larger group of diffeomorphisms that preserve the boundary $H=\{x=0\}$ (also allowing the addition of constants to $f$). The two classifications use the same objects and differ only in the group, so they can be compared directly.

This comparison is presented below.
\begin{center}
\renewcommand{\arraystretch}{1.3}
\begin{tabular}{@{}p{4.2cm}p{3.5cm}p{2.9cm}@{}}
\hline
condition on $f$ & $\Gb$-model & Arnold model\\
\hline
$f_t(0)\neq0$
 & $f_0(t)$ 
 & $t$ (nonsingular)\\[2pt]
$f_t(0)=0,f_z(0)\neq0$, $\operatorname{ord}(f_0-f_0(0))=k$
 & $f_0(t)+z$ 
 & $A_{k-1}$: $\ \pm z+t^{k}$\\
 \hline
$\dd_0f=0$, $a\neq0$, $\Delta<0$
 & $a\,t^2+\operatorname{sgn}(a)\,z^2$ 
 & $B_2$: $\ t^2+z^2$\\[2pt]
$\dd_0f=0$, $a\neq0$, $\Delta>0$
 & $a\,t^2-\operatorname{sgn}(a)\,z^2$ 
 & $B_2$: $\ t^2-z^2$\\[2pt]
$\dd_0f=0$, $a\neq0$, $\Delta=0$
 & $a\,t^2+(\text{higher terms})$ 
 & $B_\mu$: $\ t^2\pm z^\mu$, $\mu\ge3$\\
 \hline
$\dd_0f=0$, $a=0$, $b\neq0$, $\operatorname{ord}f_0=k\ge3$
 & $zt+(\text{terms in }f_0)$ 
 & $C_k$: $\ zt+t^{k}$\\[2pt]
$\dd_0f=0$, $a=b=0$, $c\neq0$, $\operatorname{ord}f_0=3$
 & $\pm cz^2+(terms\ in\ f_0)$ 
 & $F_4$: $\ z^2+t^3$\\
 \hline
\end{tabular}
\end{center}
The first column lists the conditions imposed on the germ $f$. The middle column lists, for each condition on $f$, our $\Gb$-model: on the strata $S_0,S_1$ this is the complete normal form of \ref{thm:regular}; in the critical rows ($d_0f=0$) it is only the model of the $2$-jet from \ref{prop:jets}. The last column is Arnold's model for the same $f$.

\end{remark}


\section{\texorpdfstring{Embedding of $b$-symplectomorphisms in flows of b-symplectic fields}{Embedding of b-symplectomorphisms in Hamiltonian flows}}\label{sec:embedding}

As we saw, the flow of a b-symplectic vector field consists of b-symplectomor--phisms. In this section, following \cite{Moser}, we ask ourselves if all b-symplectomor--phisms can be embedded in the flow of a b-symplectic vector field, to find that the answer is negative.

\begin{theorem}\label{thm:E}
Let $\varphi\in\Gb$ with $u_0(0)>0$ and let $F$ be a smooth function, defined in a neighbourhood of $0\in \mathbb{R}^2$, with $F\circ\varphi=F$ and $\dd_0 F\neq0$. Then there exist a smooth function $G$ of one variable and a neighbourhood $N$ of $0\in \mathbb{R}^2$, such that the flow of $X_{G\circ F}$ exists in $N$, up to time $1$, and $\varphi=\exp\bigl(X_{G\circ F}\bigr)$ on $N$. Conversely, if $\varphi=\exp(X_H)$ near $0$ for a smooth function $H$, then $H\circ\varphi=H$ and $\rho(\varphi)=e^{H_t(0,\cdot)}$.
\end{theorem}
\begin{proof}
Let us suppose first that $F$ belongs in the stratum $S_0$ and set $F_0=F|_Z$. We have that $F_0'(0)=F_t(0)\neq0$ and $F_0$ is a local diffeomorphism of $(\R,0)$.

By \ref{thm:regular}(a) there exists a $\Psi\in\Gb$ with $\tilde F:=F\circ\Psi^{-1}=F_0(t)$. If we define $\tilde\varphi:=\Psi\varphi\Psi^{-1}\in\Gb$, it is $\tilde F\circ\tilde\varphi=\tilde F$ and $\rho(\tilde\varphi)=\rho(\varphi)=u_0$.

Write $\tilde\varphi=(\tilde F_1,T)$. The relation $F_0(T)=F_0(t)$ and the injectivity of $F_0$ near $0$ give us $T(z,t)=t$ and we can then use the Jacobian equation to obtain $u_z=0$ and $u=u(0,t)=:a(t)$, where $\tilde{F}_1=zu$.
Hence $\tilde\varphi(z,t)=(a(t)z,\,t)=L_a$ with $a=u_0=\rho(\varphi)$.

The function $\tilde H(z,t):=\int_0^t\ln a(s)\,\dd s$ is well defined, since  $a(0)>0$. The b-symplectic vector field corresponding to this function is $X_{\tilde H}=\ln a(t)\,z\,\partial_z$, whose flow equals $(z,t)\mapsto(e^{s\ln a(t)}z,\,t)$, defined for all $s$. At $s=1$, this gives us the map $L_a$.

Consider now the function of one variable $G(r):=\bigl(\int_0^{r}\ln a\bigr)\circ F_0^{-1}$. It is $\tilde H=G\circ\tilde F$, so that $H:=\tilde H\circ\Psi=G\circ F$ and by \ref{lem:equiv}:
\[
\exp(X_{G\circ F})=\Psi^{-1}\circ\exp(X_{\tilde H})\circ\Psi=\Psi^{-1}L_a\Psi=\varphi.
\]
The flow exists up to time one because $X_{\tilde H}$ has a global flow.

Now, consider the case where $F$ belongs to the stratum $S_1$. By \ref{thm:regular}(b) and \ref{lem:equiv} we may assume $F=\tilde F=f_0(t)+z$ with $f_0'(0)=0$, and $\varphi=(zu,\,t+zS)$. In the coordinates $(c,t)$, $c=f_0(t)+z$, we have $z=c-f_0(t)$, $\omega=\frac1z\,\dd c\wedge\dd t$, and, since $F\circ\varphi=F$, $\varphi(c,t)=(c,g_c(t))$, where $g_c(t)=t+zS$.

Since $\varphi^*\omega=\omega$ off $Z$, $\frac1{z\circ g_c}\,g_c'=\frac1z$, i.e.:
\begin{equation}\label{eq:gprime}
g_c'(t)=\frac{z(c,g_c(t))}{z(c,t)}\qquad(z\neq0).
\end{equation}

Put $m:=\sup|f_0'|$, in a small rectangle containing the origin. Then $m$ is small because $f_0'(0)=0$. Let $z\neq0$ and $t_\theta=t+\theta zS$, $\theta\in[0,1]$. Writing $f_0(t+h)-f_0(t)=h\,K(t,h)$ with $|K|\le m$, we get:
\[
z(c,t_\theta)=z\,(1-\theta S K_\theta),\qquad K_\theta=K(t,\theta zS).
\]
After shrinking the rectangle so that $m\sup|S|<\tfrac12$, we have $1-\theta SK_\theta\in(\tfrac12,\tfrac32)$. Hence $z(c,\cdot)$ has the sign of $z$ and does not vanish on the interval between $t$ and $g_c(t)$.

For $z\neq0$ define:
\[
\mathcal T(c,t):=-\int_t^{g_c(t)}\frac{\dd t''}{z(c,t'')}
=-S\int_0^1\frac{\dd\theta}{1-\theta SK_\theta}.
\]
The right-hand side is smooth on the whole rectangle, including $Z$, because the denominator is greater than (or at least equal to) $\frac12$.
Differentiating the first expression and using \eqref{eq:gprime},
\[
-\partial_t\mathcal T=\frac{g_c'(t)}{z(c,g_c(t))}-\frac1{z(c,t)}=0 .
\]
So $\partial_t\mathcal T=0$ on $\{z\neq0\}$, hence on the whole rectangle, by continuity and density. Thus $\mathcal T=\tilde{\mathcal T}(c)$, with $\tilde{\mathcal T}(c)=\mathcal T(c,0)$ smooth.

Let $G'=\tilde{\mathcal T}$ and $H=G\circ F$. Then $X_H=-\tilde{\mathcal T}(c)\,z\,\partial_t$. Fix $x=(c,t)$, with $z\neq0$, and define:
\[
\Phi(t''):=-\int_t^{t''}\frac{\dd s}{z(c,s)}.
\]
Then, $\Phi$ is smooth, strictly monotone on the interval between $t$ and $g_c(t)$, and $\Phi(g_c(t))=\tilde{\mathcal T}(c)$. For $s\in[0,1]$ let $t_s$ be the point of this interval with $\Phi(t_s)=s\tilde{\mathcal T}(c)$. By differentiation we get:
\[
-\dot t_s/z(c,t_s)=\tilde{\mathcal T}(c),
\]
so $t_s$ solves the ODE of $X_H$ with $t_0=t$. By uniqueness it is the flow, which therefore exists up to time one, and $t_1=g_c(t)$. So $\exp(X_H)(x)=\varphi(x)$ for $z\neq0$. On $Z$ both maps are the identity. Hence $\varphi=\exp(X_{G\circ F})$ near $0$ and we are done.

To prove the converse, let $\varphi=\exp(X_H)=\varphi_1$, where $\varphi_s$ is the flow of $X_H$ near $0$. We have that $X_H(H)=zH_tH_z-zH_zH_t=0$, 
so:
\[
\frac{\dd}{\dd s}H\circ\varphi_s=(X_H H)\circ\varphi_s=0.
\]
Hence $H\circ\varphi_s=H$ for all $s$, in particular for $s=1$ (this is \ref{lem:flow}(i) with $f=H$).

\ref{lem:flow} applied to $f=H$ gives $\varphi_s\in\Gb$ and $\rho(\varphi_s)=e^{s\lambda}$ with $\lambda(t)=H_t(0,t)$.
Taking $s=1$ gives $\rho(\varphi)=e^{H_t(0,\cdot)}$, which is positive. In particular $u_0(0)>0$.
\end{proof}

We see that, if $\varphi \in \Gb$ is the time-1 map of b-symplectic vector field $X_H$, for some function $H$, then it necessarily has $u_0=\rho(\phi)>0$ and, in this case, it possesses the smooth first integral $H$, with $H_t(0,\cdot)=\log u_0(\cdot)$ on $Z$. Moreover, if $u_0(0)<0$, then $\varphi$ is not the time-1 map of a vector field tangent to $Z$.

\begin{proposition}\label{prop:embedded}
Let $\varphi=(zu,\,t+zS)\in\Gb$ be the time-one map of $X_H$ near $0$, with $H$ smooth and $H(0)=0$. Put $\lambda:=\log u_0$ and $f_0(t):=\int_0^t\lambda$. Then $H(0,t)=f_0(t)$, $H_t(0)=\lambda(0)$, and:
\[
S(0,0)=-H_z(0)\,\frac{e^{\lambda(0)}-1}{\lambda(0)}
\]
(the quotient is equal to $1$, if $\lambda(0)=0$). Consequently:
\begin{enumerate}[label=(\alph*),nosep]
\item if $\lambda(0)\neq0$, then $\varphi$ is conjugate in $\Gb$ to $L_{u_0}=\exp X_{f_0}$;
\item if $\lambda(0)=0$ and $S(0,0)\neq0$, then $\varphi$ is conjugate in $\Gb$ to $\exp X_{f_0+z}$.
\end{enumerate}
In particular, two such germs $\varphi,\varphi'$ with $(\lambda(0),S(0,0))\neq(0,0)$ and $(\lambda'(0),S'(0,0))\neq(0,0)$ are conjugate in $\Gb$ if and only if $\rho(\varphi)=\rho(\varphi')$.
\end{proposition}

\begin{proof}
By \ref{lem:flow}(iii), $u_0=e^{H_t(0,\cdot)}$, so $H_t(0,t)=\lambda(t)$ and, as $H(0)=0$, $H(0,t)=f_0(t)$. The linear part of $X_H=zH_t\partial_z-zH_z\partial_t$ at $0$ is:
\[
A=\begin{pmatrix}\lambda(0)&0\\-H_z(0)&0\end{pmatrix}\ \text{and}\  D\varphi(0)=e^{A}=\begin{pmatrix}e^{\lambda(0)}&0\\-H_z(0)\frac{e^{\lambda(0)}-1}{\lambda(0)}&1\end{pmatrix}.
\]
Comparing with:
\[
D\varphi(0)=\begin{pmatrix}u_0(0)&0\\S(0,0)&1\end{pmatrix}
\]
gives the formulas.

Now, if $\lambda(0)\neq0$, then $H_t(0)\neq0$ and \ref{thm:regular}(a) gives $\Psi\in\Gb$ with $\Psi_*X_H=X_{f_0}$; hence $\Psi\varphi\Psi^{-1}=\exp(\Psi_*X_H)=\exp X_{f_0}=L_{u_0}$ (\ref{lem:flow}). If $\lambda(0)=0$ and $S(0,0)\neq0$, then $H_t(0)=0$, $H_z(0)=-S(0,0)\neq0$, and \ref{thm:regular}(b) gives $\Psi_*X_H=X_{f_0+z}$. The last statement follows, since $f_0$ is determined by $\rho(\varphi)$ and, when $\lambda(0)=0$, both germs are conjugate to $\exp X_{f_0+z}$.
\end{proof}

We see that, among germs with $u_0>0$ and $(\lambda(0),S(0,0))\neq(0,0)$, $\rho$ is a complete invariant on the embeddable ones, and the classification reduces to the question of embeddability. The complementary case $\lambda(0)=S(0,0)=0$ corresponds to $\dd H(0)=0$, i.e.\ to the critical stratum of Section~\ref{sec:vf}.

Also, when $u_0(0)>0$, \ref{thm:A} can be deduced from the proposition above. Indeed, the function $\pi$ of \ref{prop:A} is invariant with $\dd\pi(0)=\dd t\neq0$, so $\varphi$ is embeddable by \ref{thm:E}, and \ref{prop:embedded}(a) gives $\varphi\sim L_{u_0}$. The direct construction we presented in section 3 is, nevertheless, needed when $u_0(0)<0$, since such germs are not embeddable.

\section{Conclusions}

As we saw, for germs of $b$-symplectic vector fields vanishing at the origin, the classification is reduced to that of functions $f\in\mathfrak m$ under $\Gb$. The case where these functions belong to the strata $S_0$, $S_1$ is completely solved, while outside $S_0\cup S_1$ the orbits have infinite codimension.

For $b$-symplectomorphisms, the normally hyperbolic case is also solved. Moreover $N(\varphi)=0$ for such germs, while germs which are not normally hyperbolic are of infinite codimension.

Two questions remain open. First, \ref{thm:E} shows that a germ with a first integral $F$ (with $\dd_0F\neq0$) and $u_0(0)>0$ is the time-one map of a vector field and then \ref{prop:embedded} solves the conjugacy problem for embedded germs with $(\lambda(0),S(0,0))\neq(0,0)$. Is every germ with $u_0(0)>0$ embedded in a flow (and thus its classification reduces to that of $b$-symplectic vector fields), at least in the remaining non-degenerate situations?

Also, is it possible to extend the results obtained here to $\R^{2n},n>1$? In our case, the behaviour of the dynamical systems on $Z$ is trivial and this allows us to control the behaviour in the direction normal to $Z$. In higher dimensions, however, non-trivial behaviour is expected on $Z$ and further conditions should be imposed on the systems under study.

We hope to comment on these questions in a future work.


\end{document}